\documentclass[12pt]{amsart}

\usepackage[foot]{amsaddr}

\begin{document}
\vspace*{-0.5cm}
\subjclass[2010]{32D15, 32H04, 32H40} 

\def\dbl{[\hskip -1pt[}
\def\dbr{]\hskip -1pt]}

\title {Holomorphic extension of meromorphic mappings along real analytic hypersurfaces}
\author{Ozcan Yazici}

\address{ Department  of Mathematics,  Middle East  Technical University,  06800 Ankara,  Turkey} 
\email{oyazici@metu.edu.tr}



\def\Label#1{\label{#1}}

\def\1#1{\ov{#1}}
\def\2#1{\widetilde{#1}}
\def\6#1{\mathcal{#1}}
\def\4#1{\mathbb{#1}}
\def\3#1{\widehat{#1}}
\def\7#1{\mathscr{#1}}
\def\K{{\4K}}
\def\LL{{\4L}}

\def \MM{{\4M}}

\def \B{{\4B}^{2N'-1}}

\def \H{{\4H}^{2l-1}}

\def \F{{\4H}^{2N'-1}}

\def \LL{{\4L}}

\def\Re{{\sf Re}\,}
\def\Im{{\sf Im}\,}
\def\id{{\sf id}\,}

\def\s{s}
\def\k{\kappa}
\def\ov{\overline}
\def\span{\text{\rm span}}
\def\ad{\text{\rm ad }}
\def\tr{\text{\rm tr}}
\def\xo {{x_0}}
\def\Rk{\text{\rm Rk\,}}
\def\sg{\sigma}
\def \emxy{E_{(M,M')}(X,Y)}
\def \semxy{\scrE_{(M,M')}(X,Y)}
\def \jkxy {J^k(X,Y)}
\def \gkxy {G^k(X,Y)}
\def \exy {E(X,Y)}
\def \sexy{\scrE(X,Y)}
\def \hn {holomorphically nondegenerate}
\def\hyp{hypersurface}
\def\prt#1{{\partial \over\partial #1}}
\def\det{{\text{\rm det}}}
\def\wob{{w\over B(z)}}
\def\co{\chi_1}
\def\po{p_0}
\def\fb {\bar f}
\def\gb {\bar g}
\def\Fb {\ov F}
\def\Gb {\ov G}
\def\Hb {\ov H}
\def\zb {\bar z}
\def\wb {\bar w}
\def \qb {\bar Q}
\def \t {\tau}
\def\z{\chi}
\def\w{\tau}
\def\Z{{\mathbb Z}}
\def\phi{\varphi}
\def\eps{\epsilon}

\def \T {\theta}
\def \Th {\Theta}
\def \L {\Lambda}
\def\b {\beta}
\def\a {\alpha}
\def\o {\omega}
\def\l {\lambda}

\def \im{\text{\rm Im }}
\def \re{\text{\rm Re }}
\def \Char{\text{\rm Char }}
\def \supp{\text{\rm supp }}
\def \codim{\text{\rm codim }}
\def \Ht{\text{\rm ht }}
\def \Dt{\text{\rm dt }}
\def \hO{\widehat{\mathcal O}}
\def \cl{\text{\rm cl }}
\def \bS{\mathbb S}
\def \bK{\mathbb K}
\def \bD{\mathbb D}
\def \bC{\mathbb C}
\def \bL{\mathbb L}
\def \bZ{\mathbb Z}
\def \bN{\mathbb N}
\def \scrF{\mathcal F}
\def \scrK{\mathcal K}
\def \mc #1 {\mathcal {#1}}
\def \scrM{\mathcal M}
\def \cR{\mathcal R}
\def \scrJ{\mathcal J}
\def \scrA{\mathcal A}
\def \scrO{\mathcal O}
\def \scrV{\mathcal V}
\def \scrL{\mathcal L}
\def \scrE{\mathcal E}
\def \hol{\text{\rm hol}}
\def \aut{\text{\rm aut}}
\def \Aut{\text{\rm Aut}}
\def \J{\text{\rm Jac}}
\def\jet#1#2{J^{#1}_{#2}}
\def\gp#1{G^{#1}}
\def\gpo{\gp {2k_0}_0}
\def\emmp {\scrF(M,p;M',p')}
\def\rk{\text{\rm rk\,}}
\def\Orb{\text{\rm Orb\,}}
\def\Exp{\text{\rm Exp\,}}
\def\Span{\text{\rm span\,}}
\def\d{\partial}
\def\D{\3J}
\def\pr{{\rm pr}}

\def \CZZ {\C \dbl Z,\zeta \dbr}
\def \D{\text{\rm Der}\,}
\def \Rk{\text{\rm Rk}\,}
\def \CR{\text{\rm CR}}
\def \ima{\text{\rm im}\,}
\def \I {\mathcal I}

\def \M {\mathcal M}

\newtheorem{Thm}{Theorem}[section]
\newtheorem{Cor}[Thm]{Corollary}
\newtheorem{Pro}[Thm]{Proposition}
\newtheorem{Lem}[Thm]{Lemma}

\theoremstyle{definition}\newtheorem{Def}[Thm]{Definition}

\theoremstyle{remark}
\newtheorem{Rem}[Thm]{Remark}
\newtheorem{Exa}[Thm]{Example}
\newtheorem{Exs}[Thm]{Examples}

\numberwithin{equation}{section}

\def\bl{\begin{Lem}}
\def\el{\end{Lem}}
\def\bp{\begin{Pro}}
\def\ep{\end{Pro}}
\def\bt{\begin{Thm}}
\def\et{\end{Thm}}
\def\bc{\begin{Cor}}
\def\ec{\end{Cor}}
\def\bd{\begin{Def}}
\def\ed{\end{Def}}
\def\be{\begin{Exa}}
\def\ee{\end{Exa}}
\def\bpf{\begin{proof}}
\def\epf{\end{proof}}
\def\ben{\begin{enumerate}}
\def\een{\end{enumerate}}

\newcommand{\dbar}{\bar\partial}
\newcommand{\genmat}{\lambda}
\newcommand{\polynorm}[1]{{|| #1 ||}}
\newcommand{\vnorm}[1]{\left\|  #1 \right\|}
\newcommand{\asspol}[1]{{\mathbf{#1}}}
\newcommand{\Cn}{\mathbb{C}^n}
\newcommand{\Cd}{\mathbb{C}^d}
\newcommand{\Cm}{\mathbb{C}^m}
\newcommand{\C}{\mathbb{C}}
\newcommand{\CN}{\mathbb{C}^N}
\newcommand{\CNp}{\mathbb{C}^{N^\prime}}
\newcommand{\Rd}{\mathbb{R}^d}
\newcommand{\Rn}{\mathbb{R}^n}
\newcommand{\RN}{\mathbb{R}^N}
\newcommand{\R}{\mathbb{R}}
\newcommand{\bR}{\mathbb{R}}
\newcommand{\N}{\mathbb{N}}
\newcommand{\dop}[1]{\frac{\partial}{\partial #1}}
\newcommand{\vardop}[3]{\frac{\partial^{|#3|} #1}{\partial {#2}^{#3}}}
\newcommand{\br}[1]{\langle#1 \rangle}
\newcommand{\infnorm}[1]{{\left\| #1 \right\|}_{\infty}}
\newcommand{\onenorm}[1]{{\left\| #1 \right\|}_{1}}
\newcommand{\deltanorm}[1]{{\left\| #1 \right\|}_{\Delta}}
\newcommand{\omeganorm}[1]{{\left\| #1 \right\|}_{\Omega}}
\newcommand{\nequiv}{{\equiv \!\!\!\!\!\!  / \,\,}}
\newcommand{\bk}{\mathbf{K}}
\newcommand{\p}{\prime}
\newcommand{\tV}{\mathcal{V}}
\newcommand{\poly}{\mathcal{P}}
\newcommand{\ring}{\mathcal{A}}
\newcommand{\ringk}{\ring_k}
\newcommand{\ringktwo}{\mathcal{B}_\mu}
\newcommand{\germs}{\mathcal{O}}
\newcommand{\On}{\germs_n}
\newcommand{\mcl}{\mathcal{C}}
\newcommand{\formals}{\mathcal{F}}
\newcommand{\Fn}{\formals_n}
\newcommand{\autM}{{\Aut (M,0)}}
\newcommand{\autMp}{{\Aut (M,p)}}
\newcommand{\holmaps}{\mathcal{H}}
\newcommand{\biholmaps}{\mathcal{B}}
\newcommand{\autmaps}{\mathcal{A}(\CN,0)}
\newcommand{\jetsp}[2]{ G_{#1}^{#2} }
\newcommand{\njetsp}[2]{J_{#1}^{#2} }
\newcommand{\jetm}[2]{ j_{#1}^{#2} }
\newcommand{\glnc}{\mathsf{GL_n}(\C)}
\newcommand{\glmc}{\mathsf{GL_m}(\C)}
\newcommand{\glc}{\mathsf{GL_{(m+1)n}}(\C)}
\newcommand{\glk}{\mathsf{GL_{k}}(\C)}
\newcommand{\smC}{\mathcal{C}^{\infty}}
\newcommand{\anC}{\mathcal{C}^{\omega}}
\newcommand{\kC}{\mathcal{C}^{k}}



\begin{abstract} Let $M\subset \mathbb C^n$ be a real analytic  hypersurface, $M'\subset \mathbb C^N$ $(N\geq n)$ be a  strongly pseudoconvex  real algebraic hypersurface of the special form and $F$ be a meromorphic mapping in a neighborhood of a point $p\in M$   which is holomorphic in one side of $M$. Assuming some additional conditions for the mapping $F$ on the hypersurface $M$,   we proved that $F$ has a holomorphic extension to $p$. This result may be used to show the regularity of CR mappings between real hypersurfaces of different dimensions.   \end{abstract}

\maketitle

\section{Introduction}\Label{int} 
The remarkable result of Forstreni\v{c} \cite{Fo} on the classification problem of proper holomorphic mappings between unit balls states that if $f$ is proper, holomorphic map from a ball in $\mathbb C^n$ to a ball in $\mathbb C^N$ and smooth of class $C^{N-n+1}$ on the closure then $f$ is a rational mapping.  He posed the question of the holomorphic extendibility of such a rational mapping to any boundary point.  In \cite{CS}, Cima and Suffridge proved that every such mapping extends holomorphically to a neighborhood of the closed ball. This result was extended by Chiappari \cite{Ch} by replacing the unit ball in domain with an arbitrary real analytic hypersurface in $\mathbb C^n$. 

This results are also related to regularity of CR mappings between real hypersurfaces. When the real hypersurfaces lie in the complex spaces of same dimension, CR mappings of given smoothness must be real-analytic,  (see for example \cite {BER}). In the case of real hypersurfaces of different dimensions, analyticity of CR mappings with given smoothness on the boundary was shown provided that the target is a real sphere   (see for example \cite{BHR, Mir} ).  In the proof, they first show that  the CR mappings extend meromorphically. Then using the results of Chiappari and Cima-Suffridge, this meromorphic extension defines an  analytic extension.

In this work,  we obtain a holomorphic extension result for meromorphic mappings with more general target spaces. More precisely we  prove the following theorem.

\begin{Thm}\label{main} Let $M\subset \mathbb C^n$ be a real analytic hypersurface   and  $M'\subset\mathbb C^N$ be a   strongly pseudoconvex real algebraic  hypersurface which is  locally equivalent to $\Im z_N'=p(z',\bar z')$ by a birational holomorphic change of coordinates  at a point $q\in M'$, where $(z',z_N')\in \mathbb C^N$, $N\geq n$ and $p(z',\bar z')$ is a real valued polynomial.  Let $U\subset \mathbb C^n$ be a neighborhood of a point $p\in M$ and $\Omega$ be the portion of $U$ lying on one side of $M$. If $F:U\rightarrow \mathbb C^N$ is a meromorphic mapping which maps $\Omega$ holomorphically to one side of $M'$, extends continuously on $\overline \Omega$, $F(M\cap U)\subset M' $ and $F(p)=q$, then $F$ extends holomorphically to a neighborhood of $p$.   
\end{Thm}

In the statement of Theorem \ref{main}, by $F(M\cap U)\subset M'$, we mean that $\displaystyle  \lim_{\Omega \ni z\to p} F(z)\in M'$ and  $\displaystyle F(p):=  \lim_{\Omega \ni z\to p} F(z)$   for all $p\in M\cap U$.    Note that  Theorem \ref{main} improves the result of Chiappari \cite{Ch} by replacing the sphere in the target with  a  special type of  real algebraic   hypersurface.  One can not expect to have extension for mappings with arbitrary targets. There are examples of proper rational mappings from the unit ball to  a compact set that can not be extended holomorphically through the boundary, (see \cite{CS}, \cite{IM} ).

\section {Proof of Theorem \ref{main}}

\begin{proof} For simplicity, we will take $p=(0,\dots,0).$ Since the ring of germs of holomorphic functions is a unique factorization domain, we may assume that $F=\frac{f}{g}$ where $f=(f_j)_{1\leq j \leq N}$ is a holomorphic mapping and $g$ is a holomorphic function near $0\in \mathbb C^n$ which has no common factor with $f$. If $g(0)\neq 0$, then $F$ defines a holomorphic mapping near $0$. Hence we may assume that $g(0)=0$.     

Let $M$ be given by $\psi(z,\overline{z})=0$ for some real analytic function $\psi$ near $0$ such that $\frac{\partial \psi}{\partial z_1}(0)\neq 0$. We define a non-zero holomorphic function $m(z)=\sum_{i=1}^{n}m_i z_i$ where $m_i=\frac{\partial \psi}{\partial z_i}(0)$. Since the zero sets of holomorphic functions are of measure zero, we can find a  point $a=(a_1,\dots,a_n)\neq 0$ such that  $m(a)\neq 0$, $g(a)\neq 0$, $f_j(a)\neq 0$ for all  $j=1,\dots , N$. Here we have assumed that $f_j'$-s are not identically equal to 0, otherwise we can replace the those $f_j'$-s with zeros in the rest of the proof.  Now we change the coordinates by $$z_i=a_i\zeta_1+\sum_{j=2}^{n}b_{ij}\zeta_j .$$ Since $a\neq 0$, we can choose $b_{ij}$ so that $\zeta=(\zeta_1,\dots,\zeta_n)$ gives a non-singular linear change of coordinates. In these new coordinates $\zeta=(\zeta_1,\dots ,\zeta_n)$, we have that $f_j(1,0,\dots,0)=f_j(a)\neq 0$, $g(1,0\dots,0)=g(a)\neq 0$  and 
$$\frac{\partial \psi }{\partial \zeta_1}(0)=\sum_{i=1}^n\frac{\partial \psi}{\partial z_i}(0)a_i=m(a)\neq 0.$$  

For the convenience, we will denote the new coordinates by $z$ again. Then we may assume that   $f_j(z_1,0)\not \equiv 0$, $g(z_1,0)\not\equiv 0$  and $\frac{\partial \psi}{\partial z_1}(0)\neq 0$. Hence $M$ can be defined as a graph $z_1=\rho(\overline{z_1}, \tilde z, \overline {\tilde z})$ where $\tilde z=(z_2,\dots,z_n)$ and $\rho(z_1,\lambda, \tau)$ is a holomorphic function near $0$ in $\mathbb C\times \mathbb C^{n-1}\times \mathbb C^{n-1}$. We may also assume that $\rho(\overline {z_1},0,0)=\overline{z_1}.$

By Weierstrass preparation theorem, $g$ can be written as $g(z_1,\tilde{z})=u(z_1,\tilde{z})h(z_1,\tilde{z})$ where  $u(z_1,\tilde{z})=\sum_{j=0}^{l}a_j(\tilde{z})z_1^j$ is a Weierstrass polynomial so that $a_l(\tilde{z}) \equiv 1$, $a_j(0)=0$ for $0\leq j \leq l-1$ and $h(0,0)\neq 0.$  Since $F$ is bounded as $ z=(z_1,\tilde {z})\rightarrow 0$ in $\Omega$, $f_j$ can be decomposed as $f_j(z_1,\tilde{z})=u_j(z_1,\tilde{z})h_j(z_1,\tilde{z})$ where $u_j'$-s are Weierstrass polynomials in $z_1$ of degree $k_j\geq l$ and $h_j(0,0)\neq 0$. Using division algorithm, one can find   $r_j$  of degree smaller than $l$ in $z_1$ such that $u_j(z_1,\tilde{z})=u(z_1,\tilde{z})d_j(z_1,\tilde{z})+r_j(z_1,\tilde{z}).$ Setting $D_j=d_jh_j$, $R_j=r_jh_j$, $D=(D_j)_{j=1}^N$, $R=(R_j)_{j=1}^N$, we have $f=uD+R$. Our aim is to show that $R\equiv 0$.\\

Since $M'$ is strongly pseudoconvex, by a holomorphic change of variables it can be written as $$M'=\{({z'},z'_N)\in \mathbb C^N: \Im z_N'-\sum_{j=1}^{N-1}|z_j'|^2+\phi({z'},\overline{ {z'}})=0\}$$ where $\phi\equiv 0 $ or $\phi$ is a real valued polynomial of degree bigger than $2$. If $\phi\equiv 0$  then $M'$ is locally equivalent to the  real sphere in $\mathbb C^N$ and Theorem \ref{main} follows from the main result in \cite {Ch}.

  Hence we can assume that $\phi \nequiv 0$. Let's write $\phi$ as $$\phi({z'},\overline{{z'}})=\sum_{I,J}\alpha_{I,J}{z'}^{I}{\overline{z'}}^J .$$
Since $\phi$ is real valued, $\overline {\alpha_{IJ}}=\alpha_{JI}$ and hence the highest degrees of $z'$ and $\overline {z'}$  in $\phi$ are the same, say they are equal to $d\geq 2.$   \\

 Since $F$ maps $M$ into $M'$, $\forall z\in M$ we have that
\begin{eqnarray}\label{1} \frac{f_N(z)}{g(z)}- \frac{\overline{f_N(z)}}{\overline{g(z)}}-2i\sum_{j=1}^{N-1}\frac{|f_j(z)|^2}{|g(z)|^2}-2i\phi\left(\frac{\tilde f(z)}{g(z)},\frac{\overline{\tilde f(z)}}{\overline {g(z)}}\right)\equiv 0
\end{eqnarray}
where  $\tilde{f}=(f_1,\dots,f_{N-1})$.
Multiplying both sides of the above equation by $g(z)^d{\overline{g(z)}}^d$, we obtain that 

\begin{eqnarray}\label{new0} f_N(z)g(z)^{d-1}\overline {g(z)}^d-\overline{f_N(z)}g(z)^d\overline {g(z)}^{d-1}-2i\sum_{j=1}^{N-1}| f_j(z)|^2g(z)^{d-1}\overline {g(z)}^{d-1} \\  \nonumber -2i    g(z)^d{\overline{g(z)}}^d \phi\left(\frac{\tilde f(z)}{g(z)},\frac{\overline{\tilde f(z)}}{\overline {g(z)}}\right)\equiv 0 
\end{eqnarray}

For $z=(z_1,\tilde{z})$, we set $z^*=(\rho(\overline{z_1},\tilde{z},\overline {\tilde{z}}),\tilde{z})$, $s^*(z)=\overline{s(z^*)}$ for any function $s$.
Then 

\begin{eqnarray}\label{new}  f_N(z)g(z)^{d-1}\overline {g(z^*)}^d-\overline{f_N(z^*)}g(z)^d\overline {g(z^*)}^{d-1}-2i\langle \tilde f(z),\tilde f(z^*)\rangle g(z)^{d-1}\overline {g(z^*)}^{d-1}\\ \nonumber -2i g(z)^d{\overline{g(z^*)}}^d \phi\left(\frac{\tilde f(z)}{g(z)},\frac{\overline{\tilde f(z^*)}}{\overline {g(z^*)}}\right)
\end{eqnarray}
is a  holomorphic function of $z_1$  and by  (\ref{new0}) it is  equal  to $0$ whenever $z_1=\rho(\overline{z_1},\tilde z,\overline {\tilde z})$, that is  when $z=z^*$. Here $\langle \, , \rangle$ denotes the standard inner product, that is  $\langle a,b \rangle =\sum_{i=1}^N a_i\overline {b_i}$ for $a=(a_1,\dots,a_N)$ and $b=(b_1,\dots,\ b_N)$ in $\mathbb C^N$.  For a fixed $\tilde z_0$, the real codimension of the set  $\{z_1=\rho(\overline{z_1},\tilde z_0,\overline{\tilde z_0})\}$ in $\mathbb C^n$ is at most the sum of real codimensions  of $M$ and $\{\tilde{z}=\tilde z_0\}$. Hence the real dimension of the set $\{z_1=\rho(\overline{z_1},\tilde z_0,\overline{\tilde z_0})\}$ is at least $1$.  It follows that  the function above  is identically $0$ as a function of $z_1$.  

Using the identities $f=uD+R$, $g=uh$ and $\tilde f(z^*)=\overline {\tilde f^*(z)}$, it follows from (\ref{new}) that 

\begin{eqnarray}\label{2} &&u^d{u^*}^d(D_Nh^{d-1}{h^*}^d-D_N^*h^d{h^*}^{(d-1)}-2i\langle \tilde D,\overline{\tilde D^*}\rangle h^{d-1}{h^*}^{(d-1)})\\ && +u^{d-1}{u^*}^d (R_N h^{d-1}{h^*}^d -2ih^{d-1}{h^*}^{(d-1)}\langle \tilde R,\overline{{\tilde D}^*} \rangle) \nonumber \\ &&+u^d{u^*}^{(d-1)}(R_N^*h^d{h^*}^{(d-1)}  -2ih^{d-1}{h^*}^{(d-1)}\langle \tilde D,\overline{\tilde R^*}\rangle ) \nonumber \\ &&-2iu^{d-1}{u^*}^{(d-1)}h^{d-1}{h^*}^{(d-1)}\langle \tilde R,\overline {\tilde R^*}\rangle-2ig(z)^d{\overline{g(z^*)}}^d \phi\left(\frac{\tilde f(z)}{g(z)},\frac{\overline{\tilde f(z^*)}}{\overline {g(z^*)}}\right)\equiv 0. \nonumber
\end{eqnarray}
where $\tilde D=(D_1,\dots,D_{N-1})$ and $\tilde R=(R_1,\dots,R_{N-1})$.

Let $\tilde z=0$. We note that  $$u^*(z_1,0)=\overline{ u(\rho(\overline {z_1},0,0),0)}=\overline{\rho(\overline {z_1},0,0)^l}=z_1^{l} $$  and $u(z_1,0)=z_1^l$.  Let us assume that $R(z_1,0)\nequiv 0$ and the multiplicity of $z_1$ in $R(z_1,0)$ is $a$ for some $0\leq a < l.$ That is  $R(z_1,0)=z_1^a Q(z_1)$ for some holomorphic function $Q$ such that $Q(0)\neq 0$. The multiplicity of $z_1$ in the first summand of the function in (\ref 2) is greater than or equal to $2dl$. In the second  and the third summands, the multiplicity of $z_1$ is greater than or equal  to $(2d-1)l+a$. In the fourth summand, the multiplicity of $z_1$  is greater than or  equal to $2(d-1)l+2a$. The equation $(\ref{2})$ implies that   
\begin{eqnarray}\label{111}\min\{2dl,(2d-1)l+a,2(d-1)l+2a \}=2(d-1)l+2a 
\end{eqnarray}
 must be smaller than or equal to the multiplicity of  $z_1$ in  the last term  $g(z)^d{\overline{g(z^*)}}^d \phi\left(\frac{\tilde f(z)}{g(z)},\frac{\overline{\tilde f(z^*)}}{\overline {g(z^*)}}\right)$. 

Note that the multiplicity of $z_1$ in $f=uD+R$  and  in $g=uh$ are  $a$ and $l$, respectively.  By writing $$g(z)^d{\overline{g(z^*)}}^d \phi\left(\frac{\tilde f(z)}{g(z)},\frac{\overline{\tilde f(z^*)}}{\overline {g(z^*)}}\right)=\sum_{|I|,|J|\leq d}\alpha_{IJ}\tilde f(z)^Ig(z)^{d-|I|}{\overline{\tilde f(z^*)}}^J{\overline{g(z^*)}}^{d-|J|}$$
we see that the multiplicity of $z_1$ in the last summand in $(\ref{2})$ is equal to 
\begin{eqnarray} \label{112} \min_{|I|,|J|, \alpha_{IJ} \ne 0}\{a|I|+l(d-|I|)+a|J|+l(d-|J|) \}\leq \min_{|J|}\{ad+ld+|J|(a-l) \}. 
\end{eqnarray}
The inequality above is obtained by taking $|I|=d$.
We have the following cases for $d$. \\

Case 1: $d=2$. Since the total degree of  $\phi$ is bigger than or equal to $3$, when $|I|=2$, $|J|$ must be at least one. Then it follows from $(\ref{111})$ and  $(\ref{112})$ that $$2l+2a\leq  \min_{|J|}\{2a+2l+|J|(a-l) \}\leq 3a+l.$$ The second inequality above follows from the fact that  $|J|\geq 1$ and $a-l<0$. But this implies that $l\leq a$, which contradicts to the choice of $a$ and hence $R(z_1,0)\equiv 0.$    \\

Case 2: $d>2$. It follows from  $(\ref{111})$ and  $(\ref{112})$ that  $$2(d-1)l+2a\leq   \min_{|J|}\{ad+ld+|J|(a-l) \}\leq ad+ld.$$ But this implies that $d\leq 2$. Hence again we have that $R(z_1,0)\equiv 0$.

 Now we suppose that $R\not\equiv 0$.  We may assume that $R(a)\neq 0$ for some $a=(a_1,\dots,a_n)$ such that $a_2\neq 0$. We change the coordinates from $z$ to $\zeta$ defined by 
  $$z_1=a_1\zeta_2+ \zeta_1,\; z_i=a_i\zeta_2+\sum_{j=3}^nb_{ij}\zeta_j, \; i=2,\dots,n.$$ 
  
  Since $(a_2,\dots,a_n)\neq 0$, $b_{ij}$ can be chosen so that $\zeta$ gives a non-singular linear change of coordinates. In these new coordinates, $R(0,1,0,\dots,0)=R(a)\neq 0$, $R(\zeta_1,0)\equiv 0$, $f_j(\zeta_1,0)\not \equiv 0$, $g(\zeta_1,0)\not \equiv 0$ and $\frac{\partial \phi}{\partial \zeta_1}(0)=\frac{\partial \phi}{\partial z_1}(0)\neq 0$.  We denote these new coordinates by $z$ again. Then  $$R(z_1,0)\equiv 0, R(z_1,z_2,0\dots,0)\not\equiv 0,$$ $f_j$ and $g$ do not vanish on $z_1$-axis and $M$ can be written as $z_1=\rho(\overline {z_1}, \tilde z, \overline {\tilde z} )$ near the origin. 
  
Since $R(z_1,z_2,0,\dots,0)$ is a non-zero analytic function of $z_2$  vanishing at $z_2=0$  there exists the largest  integer $k\geq 1$ such that $z_2^k$ divides $R(z_1,z_2,0,\dots,0)$. We define $G_i=\frac{R_i}{z_2^k}$, $G=(G_1,\dots,G_N)$ and $\tilde G=(G_1,\dots, G_{N-1})$.  Note that 

\begin{eqnarray} \label{113} G(z_1,0,\dots ,0)\neq 0.
\end{eqnarray}  

Then dividing the terms in (\ref{2}) by $|z_2|^{2k}$ , we obtain  that

\begin{eqnarray}\label{12}  
&& \frac{1}{|z_2|^{2k}}u^d{u^*}^d(D_Nh^{d-1}{h^*}^d-D_N^*h^d{h^*}^{(d-1)}-2i\langle \tilde D,\overline{\tilde D^*}\rangle h^{d-1}{h^*}^{(d-1)})\\  && +\frac{1}{{\overline{z_2}}^k} u^{d-1}{u^*}^d (G_N h^{d-1}{h^*}^d -2ih^{d-1}{h^*}^{(d-1)}\langle \tilde G,\overline{{\tilde D}^*} \rangle) \nonumber \\ &&+\frac{1}{{z_2}^k}u^d{u^*}^{(d-1)}(G_N^*h^d{h^*}^{(d-1)}  -2ih^{d-1}{h^*}^{(d-1)}\langle \tilde D,\overline{\tilde G^*}\rangle ) \nonumber \\ &&-2iu^{d-1}{u^*}^{(d-1)}h^{d-1}{h^*}^{(d-1)}\langle \tilde G,\overline {\tilde G^*}\rangle-\frac{2i}{|z_2|^{2k}}g(z)^d{\overline{g(z^*)}}^d \phi\left(\frac{\tilde f(z)}{g(z)},\frac{\overline{ \tilde f(z^*)}}{\overline {g(z^*)}}\right)\equiv 0. \nonumber
\end{eqnarray}

We take $(z_3,\dots,z_n)=0$ and let $z_2\rightarrow 0$ in above equation. Considering the order of $z_1$ in all terms in (\ref{12}), as in above argument for $R(z_1,0,\dots,0)$, we obtain that  $G(z_1,0\dots, 0)=0$ when  $z_1=\rho(\overline{z_1},0,0)$. But this  contradicts to $(\ref{113})$.  Consequently, $R\equiv 0$ and $F=\frac{f}{g}=\frac{D}{h}$ defines a holomorphic mapping near $0\in\mathbb C^n.$
\end{proof}

\noindent \textbf{Acknowledgments.} I am grateful to Nordine Mir for his suggestion to work on this problem and for useful discussions on this subject. I would like to thank the referee for his/her remarks and suggestions which helped to improve the presentation of the paper. The author is supported by T\"UB\.ITAK 2232 Proj. No 117C037.


\begin{thebibliography}{999}

\bibitem {BER}  Baouendi, M.S. ;   Ebenfelt, P. ;   Rothschild,  L.P. : Real submanifolds in complex space and their mappings, Princeton Math. Series {\bf 47}  Princeton Univ. Press, (1999).

\bibitem {BHR}  Baouendi, M.S.;  Huang, X. ;  Rothschild,  L.P. : Regularity of CR mappings between algebraic hypersurfaces, Invent. Math. {\bf 125} (1996), 13-36.

\bibitem  {Ch} Chiappari, S. :  Holomorphic extension of proper meromorphic mappings, {\em Michigan Math. J.}  {\bf 38} (1991), 167-174.

\bibitem {CS} Cima, J. A. ; Suffridge, T. J.  Boundary behavior of rational proper maps, {\em Duke Math. J.} {\bf 60} (1990), no. 1, 135-138.

\bibitem  {Fo} Forstneri\v{c}, F. :  Extending proper holomorphic mappings of positive codimension, {\em Invent. Math.} {\bf 95} (1989), no. 1, 31-61

\bibitem  {IM} Ivashkovich, S. ;  Meylan, F. :  An example concerning holomorphicity of meromorphic mappings along real hypersurfaces, {\em  Michigan Math. J.} {\bf 64} (2015), no. 3, 487-491.

\bibitem {Mir}  Mir, N. :  Analytic regularity of CR maps into spheres, {\em  Math. Res. Lett.} {\bf 10} (2003), no. 4, 447-457. 


\end{thebibliography}
\end{document}